\documentclass[11pt]{amsart}
\usepackage{amssymb,amsmath,amscd,amsfonts,amsthm,mathrsfs}
\usepackage{verbatim,paralist}
\usepackage{enumerate}
\RequirePackage{ifthen}
\usepackage[english]{babel}
\usepackage{xcolor}
 \usepackage{tikz}
\usepackage{latexsym}
\usepackage{amsmath}
\usepackage{amsfonts}
\usepackage{amssymb}
\usepackage{latexsym}

\usepackage{color}

\newcommand{\Ac}{\mathcal{A}}

\newcommand{\Wc}{\mathcal{W}}

\newcommand{\Oc}{\mathcal{O}}
\newcommand{\Sc}{\mathcal{S}}
\newcommand{\Vc}{\mathcal{V}}

\newtheorem{theorem}{Theorem}[section]
\newtheorem{proposition}[theorem]{Proposition}
\newtheorem{lemma}[theorem]{Lemma}
\newtheorem{remark}[theorem]{Remark}
\newtheorem{example}[theorem]{Example}
\newtheorem{definition}[theorem]{Definition}
\newtheorem{corollary}[theorem]{Corollary}

\title{Fredholm theory in quaternionic Banach algebras}
\author{Hatem Baloudi}
\address{Hatem Baloudi, Department of Mathematics, Faculty of Sciences of Gafsa, University of Gafsa, 2112 Zarroug,
Tunisia}

\subjclass[2010]{46S10, 47A60, 47A10, 47A53, 47B07}
\keywords{Quaternions, Quaternionic Banach algebra, Weyl spectrum, Fredholm spectrum.}

\begin{document}

\maketitle
\begin{abstract}Muraleetharan and Thirulogasanthan in (J. Math. phys. {\bf 59}, No. 10, 103506, 27p. (2018)) introduced the concept of Calkin S-spectrum of a bounded quaternionic linear operators. The study of this spectrum is establisched using the Fredholm operators theory. Motivated by this, we study the general framework of the Fredholm element with respect to a quaternionic Banach algebra homomorphism. First, we investigate the Fredholm S-spectrum of the sum of two elements in quaternionic Banach algebra by means of the Fredholm S-spectrum of the two elements. Next, we prove a perturbation result on this spectrum. We also study the boundary S-spectrum. As application, we investigate the Fredholm and Weyl S-spectra of bounded right quaternionic linear operators.\end{abstract}

\tableofcontents
\section{Introduction}
In the complex setting, much attention has been paid to Fredholm theory of  linear operators on a complex Banach spaces \cite{BJ,AN,RH,J1}. It is useful for the study of non-parabolicity at infinity of Dirac type operators on non-compact Riemannian manifuld \cite{Carron}. We refer to \cite{Ayedi} for the study of this point on the discrete Gauss-Bonnet operator on an infinite graph. Moreover, Fredholm operators theory plays an important role in the study of non-discrete spectrum, and  the study of the stability of some parts on the spectrum under some perturbation \cite{HB,J2, Kato,Schechter}.
\vskip 0.1 cm

In complex spectral theory, it is well known (Atkinson's theorem) that a bounded linear operator on a complex Banach space is Fredholm if, and only if, it is invertible in the Calkin algebra. In particular, the Fredholm spectrum of a bounded  operator $A$ is given by the spectrum of $\pi(A)$, where $\pi$ is the natural quotient map associated with the Calkin algebra. As an extension, here Fredholm and perturbation theory relative to a unital homomorphism $\Ac:\ \Vc\longrightarrow\Wc$ between unital complex Banach algebras $\Vc$ and $\Wc$. This class  introduced firstly by R. Harte in \cite{RH}. A generalized Fredholm theory is recently discuss in \cite{3}.
\vskip 0.1 cm

In the quaternionic setting, the quaternionic  multiplication is not commutative. This leads to three types of Hilbert spaces: left, right, and two-sided, depending on how vectors are multiplied by scalars. Therefore, an apparent problem in defining the spectrum of a quaternionic operator. In fact, if $T$ is a right linear operator on a quaternionic vector space $\Vc$, the function $x\longmapsto Tx-xq$ is not right linear, we refer to \cite{NCFCBO, FJDP} for this point of view. The fundamental suggestion of \cite{NCFCBO} is to define the spectrum using the Cauchy kernel series $S^{-1}(q,T)=\sum_{n}T^{n}q^{-n-1}$, where $q$ is a quaternionic number and $T$ is a right linear operator on a right quaternionic vector space. More precisely, they identified the operator whose inverse is $S^{-1}(q,T)$. This leads to define a new spectrum called S-spectrum. Based on this concept, \cite{BK} introduced and studied the Fredholm and Weyl S-spectra.
\vskip 0.1 cm

The aim of this article is to complete the knowledge of the concept of the Fredholm S-spectrum and give a general framework of Fredholm theory relative to a quaternionic Banach algebra homomorphism. In general, the set of right linear operators on the quaternionic Banach space is not quaternionic Banach algebra with respect to the composition of operators. If $\Vc$ is a right separable Hilbert space, it is possible to define the left scalar multiplication on $\Vc$ by using an arbitrary Hilbert basis on $\Vc$. We refer to \cite{RVA} for an explanation of this construction. Thanks to this construction, the set of right linear operators becomes a quaternionic two-sided Banach $C^{*}-$algebra with unity, see \cite[Theorem 3.4]{RVA}. In this regard, \cite{BK} introduce the Calkin spherical spectrum as the spherical spectrum of the quotient map image of a bounded right linear operator on the Calkin algebra. They gave a characterization of this new spectrum by using the concept of Fredholm operators. Motivated by this, we introduce and we study the general framework of the Fredholm element with respect to a quaternionic Banach algebra homomorphism. In particular, we give some new results for the Calkin and Weyl S-spectra of the bounded right quaternionic linear operator.
\vskip 0.1 cm

To start off, we introduce the quaternionic version of the Fredholm elements relative to the homomorphism $\Ac:\ \Vc\longrightarrow\Wc$, where $\Vc$ and $\Wc$ are two quaternionic two-sided Banach algebra with unit $1\neq 0$. In addition, we define the Fredholm and Weyl S-spectra with respect to $\Ac$. One of the main results is: The Fredholm S-spectrum of the sum of two elements in a Banach algebra is written by means of the Fredholm S-spectra of the two elements when their products are in the kernel of $\Ac$, see Theorem \ref{sum}. This result is applicable for the Calkin (or Fredholm) S-spectrum of the sum of two bounded quaternionic right linear operators when their products are compact, see Theorem \ref{sum1}. On the other hand, we show the invariance of the Fredholm and Weyl S-spectra under certain perturbations. Furthermore, we investigate the boundary S-spectrum of an element in quaternionic Banach algebra.
\vskip 0.1 cm

The article is organized as follows. In Section 2, we gather some definitions and notations, the facts connected to our work. In Section 3, the concept of the Fredholm and Weyl  S-spectra in a quaternionc Banach algebra is introduced. Some properties and perturbation results of these spectra are investigated. In Section 4, we introduce and we study the boundary  S-spectrum. Finally, in Section 5, we apply the results of Section 3 to describe the Fredholm and Weyl S-spectra of the sum of two bounded right linear operators. We also discuss the coincidence between the Fredholm and Weyl S-spectrum.

\section{Mathematical preliminaries}
In this section we recall some standard definitions and we gives some results that we will need in the sequel. For details we refer to the reader to \cite{SL,FJDP,FIDC1,RVA}.
\subsection{Quaternions}
The algebra of quaternion $\mathbb{H}$ is defined as
\begin{align*}\mathbb{H}=\Big\{q=q_{0}+iq_{1}+jq_{2}+kq_{3}:(q_{0};q_{1},q_{2},q_{3})\in\mathbb{R}^{4}\Big\},\end{align*}
where the imaginary units $i,\ j,\ k$ satisfy
\begin{align*}i^{2}=j^{2}=k^{2}=-1,\ ij=-ji=k,\ jk=-kj=i,\ ki=-ik=j.\end{align*}
The quaternionic conjugate of $q=q_{0}+iq_{1}+jq_{2}+kq_{3}$ is
\begin{align*}\overline{q}=q_{0}-iq_{1}-jq_{2}-kq_{3},\end{align*}
while $|q|=(q\overline{q})^{\frac{1}{2}}$ denotes the usual norm of the quaternion $q$. If $q$ is non-zero element, it has inverse
\begin{align*}\displaystyle q^{-1}=\frac{\overline{q}}{|q|^{2}}.\end{align*}
A purely imaginary quaternion with absolute value $1$ is called an imaginary unit. We denote the set of all imaginary units by $\mathbb{S}$, that is,
\begin{align*}\mathbb{S}=\left\{q\in\mathbb{H}:\ \overline{q}=-q,\ |q|=1\right\}.\end{align*}
The name imaginary unit is justified by the fact that, for any $I\in\mathbb{S}$, we have
\begin{align*}I^{2}=-\overline{I}I=-|I|^{2}=-1.\end{align*}
We set
\begin{align*}\mathbb{C}_{I}=\Big\{x+Iy\in\mathbb{H}:\ x,y\in\mathbb{H}\Big\},\end{align*}
i.e, $\mathbb{C}_{I}$ is the real vector subspace of $\mathbb{H}$ generated by $1$ and $I\in\mathbb{S}$.
Therefore, the plane $\mathbb{C}_{I}$ is isomorphic to the field of complex number $\mathbb{C}$.
For a quaternion $q\in\mathbb{H}$, we set
\begin{align*}I_{q}:=\left\{ \begin{array}{rl}
\noindent\frac{q}{|Im(q)|}, & \mbox{ if } Im(q)\neq 0,
\\~~
\\
 \displaystyle
\mbox{any } I\in\mathbb{S},&
  \mbox{ if } Im(q)=0.
\end{array} \right.\end{align*}
Then $I_{q}\in\mathbb{S}$ and $q\in\mathbb{C}_{I_{q}}$. More precisely, $q=Re(q)+I_{q}|Im(q)|$, where $Re(q)$ is the real part of the quaternion $q$ and $Im(q)$ is the vector part of $q$. The set
\begin{align*}[q]:=\{Re(q)+I|Im(q)|:\ I\in\mathbb{S}\},\end{align*}
is a 2-sphere of radius $|Im(q)|$ centred at the real point $Re(q)$.
\subsection{Slice hyperholomorphic functions}
In this section, we recall the concept of slice hyperholomorphic function and refer to \cite{GJ,FIF,FIDC,FIDC1} for surveys on
the matter.
\begin{definition}
A set $U\subset\mathbb{H}$ is called\\
\noindent $(i)$ axially symmetric if $[x]\subset U$ for any $x\in U$ and\\
\noindent $(ii)$ a slice domain if $U$ is open, $U\cap\mathbb{R}\neq\emptyset$ and $U\cap\mathbb{C}_{I}$ is a domain for any $I\in\mathbb{S}$.

\end{definition}
Let $I$ be an imaginary unit. We denote the Wirtinger derivatives with respect to the complex and the conjugate variable on the plane $\mathbb{C}_{I}$  by $\partial_{I}$ and $\overline{\partial_{I}}$, that is, the operators
\begin{align*}\partial_{I}=\frac{1}{2}\left(\frac{\partial}{\partial_{q_{0}}}-I\frac{\partial}{\partial_{q_{1}}}\right)
\mbox{   and   }\overline{\partial}_{I}=\frac{1}{2}\left(\frac{\partial}{\partial_{q_{0}}}+I\frac{\partial}{\partial_{q_{1}}}\right)\end{align*}
\noindent Let $U\subset\mathbb{C}_{I}$ be an open set and let $f:\ U\longrightarrow\mathbb{H}$ be a real differentiable function. If $\partial_{I}$ act on the right, then
\begin{align*}\overline{\partial}_{I}f(q)=\displaystyle\frac{1}{2}\left(\frac{\partial}
{\partial_{q_{0}}}f_{I}(q_{0}+Iq_{1})+\frac{\partial}{\partial_{q_{1}}}f_{I}(q_{0}+Iq_{1})I\right)\end{align*}
for $q=q_{0}+Iq_{1}\in U$, where $f_{I}=f|_{U\cap\mathbb{C}_{I}}$
\begin{proposition}\cite{FJDP} For any non-real quaternion $q\in\mathbb{H}\setminus\mathbb{R}$, there exist, and are unique, $q_{0},q_{1}\in\mathbb{R}$ with $q_{1}>0$, and $I\in\mathbb{S}$ such that $q=q_{0}+q_{1}I$.

\end{proposition}
\begin{definition}$($Slice-regular function \cite{FJDP,GC,C}$)$ Let $U$ be a domain in $\mathbb{H}$ and let $\Vc_{\mathbb{H}}^{R}$ be a vector space under right multiplication by quaternions. A real differentiable operator-valued function $f: U\longrightarrow \Vc_{\mathbb{H}}^{R}$ is said to be slice right regular if, for every quaternion $I\in\mathbb{S}$, the restriction of $f$ to the complex plane $\mathbb{C}_{I}=\mathbb{R}+I\mathbb{R}$ passing through the origin, and containing $1$ and $I$, has continuous partial derivatives (with respect to $q_{0}$ and $q_{1}$, every element in $\mathbb{C}_{I}$ being uniquely expressible as $q_{0}+q_{1}$I) and satisfies
\begin{align*}\overline{\partial}f(q_{0}+q_{1}I)=0\end{align*}\end{definition}
\noindent where $f_{I}=f|_{U\cap \mathbb{C}_{I}}$.
\subsection{Quaternionic two-sided Banach algebra}
We recall the notion of quaternionic two-sided Banach algebra and refer to \cite{GR,RVA,TA} for surveys on the matter. Let $\Vc$ be a quaternionic two-sided vector space. We say that $\Vc$ is a (associative) quaternionic two-sided Banach algebra if it is endowed with an associative product $\Vc\times\Vc\longrightarrow\Vc:\ (x,y)\longmapsto xy$ such that
\begin{enumerate}
\item $x(y+z)=xy+xz$ $\forall x,y,z\in\Vc$,
\item $(x+y)z=xz+yz$ $\forall x, y, z\in\Vc$,
\item $q(xy)=(qx)y$ $\forall x, y\in \Vc$, $\forall q\in\mathbb{H}$
\item $(xy)q=x(yq)$ $\forall x, y\in \Vc$, $\forall q\in\mathbb{H}$
\end{enumerate}
We say that $\Vc$ is a normed two-sided algebra provide $\Vc$ is endowed with an $\mathbb{H}-$norm such that $\|xy\|\leq\|x\|\|y\|$ for every $x,y\in\Vc$. If $\Vc$ is complete we say that $\Vc$ is a quaternionic two-sided Banach algebra. If in addition $\Vc$ is non-trivial and $\|1_{\Vc}\|=1$, then $\Vc$ is called a quaternionic two-sided Banach algebra with unit. \\
A map $\Ac:\Vc\longrightarrow\Wc$ between quaternionic two-sided Banach algebras with unity is called homomorphism if
\begin{enumerate}
\item $\Ac(u+v)=\Ac(u)+\Ac(v)$ and $\Ac(uv)=\Ac(u)\Ac(v)$,
\item $\Ac(qu)=q\Ac(u)$ and $\Ac(uq)=\Ac(u)q$,
\item $\Ac(1_{\Vc})=1_{\Wc}$
\end{enumerate}
for every $u,v\in\Vc$ and $q\in\mathbb{H}$. A mapping $v\longrightarrow v^{*}$ $(v\in \Vc)$  is called involution if it possesses the following properties:
\begin{enumerate}
\item $(u+v)^{*}=u^{*}+v^{*}$,
\item $(uv)^{*}=v^{*}u^{*}$,
\item $(u^{*})^{*}=u$,
\item $(puq)^{*}=\overline{q}u^{*}\overline{p}$,
\end{enumerate}
for every $u,v\in\Vc$ and $p,q\in\mathbb{H}$. Finally, a quaternionic two-sided Banach algebra is called a quaternionic two-sided $C^{*}-$algebra if $\|u^{*}u\|=\|u\|^{2}$ for all $u\in\Vc$

\begin{remark}\cite[Remark.3.3]{RVA}\\
\noindent $1)$ The concept of real two-sided Banach unital algebra is equivalent to the usual one of real Banach unital algebras.\\
\noindent $2)$ Let $I\in\mathbb{S}$. The concept of complex Banach unital algebras is equivalent to the notion of two-sided $\mathbb{C}_{I}-$Banach unital algebras $\Vc$ having the following property
\begin{align*}su=us \mbox{ for every }u\in\Vc \mbox{ and }s\in\mathbb{C}_{I}.\end{align*}
\end{remark}
Let $\Ac:\Vc\longrightarrow\Wc$ be a homomorphism of Banach algebra, where $\Vc$ and $\Wc$ are two quaternionic two-sided Banach algebra with unit $1\neq 0$. We use $\Vc^{-1}$ and $\Wc^{-1}$ for the groups of invertible elements in $\Vc$ and $\Wc$, it follows that $\Ac(\Vc^{-1})\subset\Wc^{-1}$. It is easy to see that $\Vc^{-1}$ is an open set. We consider the following definition of Fredholm element relative to the homomorphism $\Ac$. For the classical definition in complexe case we refer to \cite{HB,RH,TH}.
\begin{definition}
\begin{enumerate}
\item An element $v\in\Vc$ is called a Fredholm element relative to the homomorphism $\Ac$, if $\Ac(v)\subset\Wc^{-1}$.
\item An element $v\in\Vc$ is called a Weyl element relative to $\Ac$, if $v\in \Vc^{-1}+\Ac^{-1}(0)$, the sum of an invertible element and one whose image is zero.
\end{enumerate}
\end{definition}
\noindent We denote the set of Fredholm elements relative to the homomorphism $\Ac$ buy $\Phi_{\Ac}$. Let $\Phi_{\Ac}^{0}$ the set of Weyl elements relative to $\Ac$.
\begin{remark}\label{r1}
\begin{enumerate}
\item If $\Ac$ is continuous, then $\Phi_{\Ac}$ is an open set in $\Vc$. Indeed, $\Phi_{\Ac}=\Ac^{-1}(\Wc^{-1})$, where $\Wc^{-1}$ is an open set in $\Wc$.
\item It is clear that
\begin{align*}\Vc^{-1}\subset \Phi_{\Ac} \subset \Phi_{\Ac}^{0}.\end{align*}
\item Since $\Ac$ is a homomorphism, then $\Ac(xy)=\Ac(x)\Ac(y)$, thus, for all $x,y\in \Phi_{\Ac}$ the product $xy\in \Phi_{\Ac}$.
\item If $ab\in\Phi_{\Ac}$ and $ba\in\Phi_{\Ac}$, then $a\in\Phi_{a}$ and $b\in\Phi_{\Ac}$. Indeed, let $i,j\in\Wc$ such that $i\Ac(ab)=\Ac(ba)j=1$. Then $i\Ac(a)=i\Ac(a)\Ac(b)\Ac(a)j=\Ac(a)j.$
\end{enumerate}
\end{remark}
\subsection{S-spectrum}
Colombo et all in \cite{FIDC} introduced the concept of S-spectrum of a bounded quaternionic operator. This notion is related to the Cauchy kernel series in the quaternionic setting. In this subsection, we introduce the general framework of the S-spectrum of an element in a quaternionic Banach algebra.

\begin{definition}Let $v\in\Vc$ and $q\in\mathbb{H}$. The left Cauchy kernel series is given by
\begin{align*}\sum_{n\in\mathbb{N}}v^{n}q^{-1-n}\end{align*}
for $\|v\|<|q|$.
\end{definition}
By repeating the same proof of \cite[Theorem 4.2]{NCFCBO} we have: If $v-\overline{q}1_{\Vc}\in\Vc^{-1}$, then
\begin{align*}-(v-\overline{q}1_{\Vc})^{-1}(v^{2}-2Re(q)v+|q|^{2}1_{\Vc})\end{align*}
is the inverse of $\sum_{n\in\mathbb{N}}v^{n}q^{-1-n}$. Moreover, we have
\begin{align*}\sum_{n\in\mathbb{N}}v^{n}q^{-1-n}=-(v^{2}-2Re(q)v+|q|^{2}1_{\Vc})^{-1}(v-\overline{q}1_{\Vc}).\end{align*}
Therefore, it is natural to define the S-spectrum:
\begin{align*}\sigma_{S}(v)=\{q\in\mathbb{H}:\ v^{2}-2Re(q)v+|q|^{2}1_{\Vc}\not\in\Vc^{-1}\}.\end{align*}
\begin{definition}
Let $v\in\Vc$. For $q\in\mathbb{H}\setminus \sigma_{S}(v)$, we define the left S-resolvent as
\begin{align*}S_{L}^{-1}(q,v):=-(v^{2}-2Re(q)v+|q|^{2}1_{\Vc})^{-1}(v-\overline{q}1_{\Vc})\end{align*}
\end{definition}
\begin{remark} Let $v\in\Vc$ and  $q\in\mathbb{H}\setminus \sigma_{S}(v)$.\\
\begin{enumerate}
\item If $qv=vq$, then the left S-resolvent coincide with the classical resolvent; i.e,
\begin{align*}S_{L}^{-1}(q,v)=(v-q1_{\Vc})^{-1}\end{align*}
\item In analogy to \cite[Proposition 6.1.27]{DFIS}, $S_{L}^{-1}(q,v)$ is a $\Vc-$valued right slice regular function of the variable $q\in \mathbb{H}\setminus \sigma_{S}(v)$, that is, $\overline{\partial}S_{L}^{-1}(q_{0}+Iq_{1},v)=0$ for all $I\in\mathbb{S}$.
\end{enumerate}
\end{remark}
\begin{example}
We note that the set of matrices of type $n\times n$ with quaternion entries is quaternionic two-sided Banach algebra with the usual addition of matrices, the standard product of matrices, the left and the right multiplication $Aq=(a_{i,j}q)$ and $qA=(qa_{i,j})$.
\begin{enumerate}
\item \cite[Example 4.15, 1)]{NCFCBO} Consider the following $2\times 2$ quaternionic matrix
\begin{align*}A=\left(
 \begin{array}{cc}
 i & 0 \\
 0 & j \\
\end{array}
\right).\end{align*}
\noindent Then,
\begin{align*}A^{2}-2Re(s)A+|s|^{2}1_{M_{2}(\mathbb{H})}=\left(
 \begin{array}{cc}
 -1-2iRe(s)+|s|^{2} & 0 \\
 0 &  -1-2jRe(s)+|s|^{2}\\
\end{array}
\right)\end{align*}
and therefore, $\sigma_{S}(A)=\mathbb{S}$.
\item We now consider the following $3\times 3$ quaternionic matrix
\begin{align*}B=\left(
 \begin{array}{ccc}
0 & 0 & i\\
0 & j & 0\\
k & 0 & 0\\
\end{array}
\right)\end{align*}
Then $s\in\sigma_{S}(B)$ if and only if
\begin{align*}\left(
 \begin{array}{ccc}
-j+|s|^{2} & 0 & -2Re(s)i\\
0 & -1-2Re(s)j+|s|^{2} & 0\\
-2Re(s)k & 0 & j+|s|^{2}\\
\end{array}
\right)\left(
 \begin{array}{ccc}
q_{1} \\
q_{2} \\
q_{3} \\
\end{array}
\right)=\left(
 \begin{array}{ccc}
0 \\
0 \\
0\\
\end{array}
\right)\end{align*}
for some $q=(q_{1},q_{2},q_{3})\in\mathbb{H}^{3}\backslash\{0\}$. We observe that $q_{1}=0$ if and only if $q_{2}=0$, then we have
\begin{equation}\label{e4}1+2Re(s)j-|s|^{2}=0\end{equation}
or
\begin{align*}
(S)\displaystyle\left\{\begin{array}{ll}
(-j+|s|^{2})q_{1}-2Re(s)iq_{3}=0
\\
\\
-2Re(s)kq_{1}+(j+|s|^{2})q_{3}=0
\end{array}\right.
\end{align*}
The equation Eq.(\ref{e4}) gives
\begin{align*}1+2Re(s)j-|s|^{2}=0\mbox{ if and only if}\  s\ \in\mathbb{S}.\end{align*}
We now consider the system $(S)$: we have
\begin{align*}q_{1}=2Re(s)\displaystyle\frac{|s|^{2}+j}{|s|^{4}+1}iq_{3}\end{align*}
and replacing it in the second equation of $(S)$ we get
\begin{align*}\displaystyle\frac{4Re^{2}(s)}{|s|^{4}+1}(1+j|s|^{2})-(|s|^{2}+j)=0\end{align*}
which gives
\begin{align*}|s|=1\mbox{ and } Re(s)\in\{\pm\frac{\sqrt{2}}{2}\}.\end{align*}
We conclude that
\begin{align*}\sigma_{S}(B)=\mathbb{S}\cup\Big\{q\in\mathbb{H}:\ |q|=1\mbox{ and }Re(s)\in\{\pm\frac{\sqrt{2}}{2}\}\Big\}.\end{align*}
\end{enumerate}
\end{example}
\noindent {\bf Notation} We denote by $\mathbb{N}$ the non-negative integer. In particular, $0\in\mathbb{N}$.
\section{Fredholm S-spectrum}\label{s:sum}
\noindent The Fredholm S-spectrum of $v\in\Vc$ relative to the homomorphism $\Ac$, is given by
\begin{align*}\sigma_{\Sc,\Ac}^{\Phi}(v):=\sigma_{S}(\Ac(v))=\{q\in\mathbb{H}:\ R_{q}(v)\not\in \Phi_{\Ac}\}.\end{align*}
We define the Weyl S-spectrum of $v\in\Vc$ relative to the homomorphism $\Ac$ as
\begin{align*}\sigma_{\Sc,\Ac}^{\Phi^{0}}(v):=\{q\in\mathbb{H}:\ R_{q}(v)\not\in \Phi_{\Ac}^{0}\}.\end{align*}
\noindent \begin{remark}Let $v\in\Vc$.
\begin{enumerate}
\item If $q\in \sigma_{\Sc,\Ac}^{\Phi}(v)$ $($resp. $\sigma_{\Sc,\Ac}^{\Phi^{0}}(v))$ then $[q]\subset \sigma_{\Sc,\Ac}^{\Phi}(v)$ $($resp. $\sigma_{\Sc,\Ac}^{\Phi^{0}}(v))$. Therefore,
 $\sigma_{S,\Ac}^{\Phi}(v)$ and $\sigma_{S,\Ac}^{\Phi^{0}}(v)$ are axially symmetric.
\item In general, we have the following inclusion
\begin{align*}
\sigma_{S,\Ac}^{\Phi}(v)\subset\sigma_{S,\Ac}^{\Phi^{0}}(v)\subset\sigma_{S}(v).
\end{align*}
\end{enumerate}
\end{remark}
\begin{definition}
An element $p\in\Vc$ is called Fredholm perturbation if $p+t\in\Phi_{\Ac}$ for all $t\in\Phi_{\Ac}$. We denote this set by $Pr(\Phi_{\Ac})$.
\end{definition}
\begin{lemma}
$Pr(\Phi_{\Ac})$ is closed two-sided ideal of $\Ac$.
\end{lemma}
\proof The proof is the same as the proof in a complex Banach space, see\cite{AN,HB}.\qed
\begin{remark}Let $\Ac:\Vc\longrightarrow\Wc$ be a homomorphism of Banach algebra.
\begin{enumerate}
\item If $p\in Pr(\Phi_{\Ac})$, then
\begin{align*}\sigma_{\Sc,\Ac}^{\Phi}(p)\backslash\{0\}=\emptyset.\end{align*}
Indeed, let $q\in\mathbb{H}^{*}$, then $|q|^{2}1_{\Vc}\in\Phi_{\Ac}$. Since $Pr(\Phi_{\Ac})$ is closed two-sided ideal of $\Vc$, then $R_{q}(p)\in\Phi_{\Ac}$.
\item If $p\in Pr(\Phi_{\Ac})$, then $p\not\in\Vc^{-1}$. Indeed, $1_{\Vc}\not\in Pr(\Phi_{\Ac})$.
\item In general, we have the following inclusion
\begin{align*}\Ac^{-1}(0)\subset Pr(\Phi_{\Ac}).\end{align*}
\item Let $v\in\Vc$ and $q\in \mathbb{H}$ with $|q|>\|\Ac(v)\|$. Then
\begin{align*}v^{2}-2Re(q)v+|q|^{2}1_{\Vc}\in\Phi_{\Ac}.\end{align*}
Indeed, for each $n\in\mathbb{N}$, set
\begin{align*}a_{n}=|q|^{-2n-2}\sum_{h=0}^{n}q^{h}\overline{q}^{n-h}.\end{align*}
According to the proof of \cite[Theorem 4.3(a)]{RVA}, $a_{0}=|q|^{-2}$, $-2Re(q)a_{0}+|q|^{2}a_{1}=0$ and $a_{n-2}-2Re(q)a_{n-1}+|q|^{2}a_{n}=0$ for every $n>1$. Moreover, we have
\begin{align*}\sum_{n=0}^{\infty}\|\Ac(v^{n})\||a_{n}|\leq\sum_{n=0}^{\infty}\|\Ac(v)\|^{n}(n+1)|q|^{-n-2}.\end{align*}
Hence, $\sum_{n\in\mathbb{N}}\Ac(v^{n})$ converge absolutely. In this way, we have
\begin{align*}\Ac(R_{q}(v))(\sum_{n\in\mathbb{N}}\Ac(v^{n})a_{n})=(\sum_{n\in\mathbb{N}}\Ac(v^{n})a_{n})\Ac(R_{q}(v))=1_{\Wc}.\end{align*}
\end{enumerate}
\end{remark}
We first prove the following result:
\begin{theorem}\label{sum}
Let $\Ac:\Vc\longrightarrow\Wc$ be a homomorphism of Banach algebra, where $\Vc$ and $\Wc$ are two quaternionic two-sided Banach algebra with unit $1\neq 0$. Take $a,b\in\Vc$.\\
\noindent 1) If $ab\in Pr(\Phi_{\Ac})$ and $ba\in Pr(\Phi_{\Ac})$, then
\begin{equation}\label{e0}\sigma_{\Sc,\Ac}^{\Phi}(a+b)\backslash\{0\}=\left[\sigma_{\Sc,\Ac}^{\Phi}(a)\cup\sigma_{\Sc,\Ac}^{\Phi}(b)\right]\backslash\{0\}.\end{equation}
\noindent 2) If $ab\in \Ac^{-1}(0)$ and $ba\in \Ac^{-1}(0)$, then
\begin{equation}\label{e88}\sigma_{\Sc,\Ac}^{\Phi^{0}}(a+b)\backslash\{0\}\subset\left[\sigma_{\Sc,\Ac}^{\Phi^{0}}(a)\cup\sigma_{\Sc,\Ac}^{\Phi^{0}}(b)\right]\backslash\{0\}.\end{equation}
\noindent If, further, $\sigma_{\Sc,\Ac}^{\Phi^{0}}(a)=\sigma_{\Sc,\Ac}^{\Phi}(a)$, then
\begin{equation}\label{e99}\sigma_{\Sc,\Ac}^{\Phi^{0}}(a+b)\backslash\{0\}=\left[\sigma_{\Sc,\Ac}^{\Phi^{0}}(a)\cup\sigma_{\Sc,\Ac}^{\Phi^{0}}(b)\right]\backslash\{0\}.\end{equation}
\end{theorem}
\proof For $q\in\mathbb{H}^{*}$, we can write
\begin{equation}\label{e1}  \displaystyle\frac{R_{q}(b)R_{q}(a)}{|q|^{2}}=R_{q}(a+b)-ab-ba+\displaystyle\frac{b^{2}a^{2}-2Re(q)(b^{2}a+ba^{2}-2Re(q)ba)}{|q|^{2}}\end{equation}
and
\begin{equation}\label{e2}  \displaystyle\frac{R_{q}(a)R_{q}(b)}{|q|^{2}}=R_{q}(a+b)-ab-ba+\displaystyle\frac{a^{2}b^{2}-2Re(q)(a^{2}b+ab^{2}-2Re(q)ab)}{|q|^{2}}.\end{equation}
Let $q\not\in\sigma_{\Sc,\Ac}^{\Phi}(a)\cup\sigma_{\Sc,\Ac}^{\Phi}(b)\cup\{0\}$, then
\begin{align*}R_{q}(a)\in\Phi_{\Ac}\mbox{ and }R_{q}(b)\in\Phi_{\Ac}.\end{align*}
Using Remark \ref{r1} $(3)$, we obtain $R_{q}(a)R_{q}(b)\in\Phi_{\Ac}$. Since $ab$ and $ba$ are two Fredholm perturbation, we can apply Eq. (\ref{e2}), we infer that
\begin{align*}|q|^{2}R_{q}(a+b)\in\Phi_{\Ac}.\end{align*}
Also, since $q\neq 0,$ then $R_{q}(a+b)\in\Phi_{\Ac}$. Therefore,
\begin{equation}\label{e3}\sigma_{\Sc,\Ac}^{\Phi}(a+b)\backslash\{0\}\subset\left[\sigma_{\Sc,\Ac}^{\Phi}(a)\cup\sigma_{\Sc,\Ac}^{\Phi}(b)\right]\backslash\{0\}.\end{equation}
\noindent To prove the inverse inclusion of Eq. (\ref{e3}). Suppose that $q\not\in\sigma_{\Sc,\Ac}^{\Phi}(a+b)\cup\{0\}$, then $R_{q}(a+b)\in\Phi_{\Ac}$. Since $ab\in Pr(\Phi_{\Ac})$ and $ba\in Pr(\Phi_{\Ac})$, then by Eqs. (\ref{e1}) and (\ref{e2}), we have
\begin{align*}R_{q}(a)R_{q}(b)\in\Phi_{\Ac}\mbox{  and  }R_{q}(b)R_{q}(a)\in\Phi_{\Ac}.\end{align*}
Again, using Remark \ref{r1} $(4)$, we have $R_{q}(a)\in\Phi_{\Ac}$ and $R_{q}(b)\in\Phi_{\Ac}$. Hence $q\not\in \sigma_{\Sc,\Ac}^{\Phi}(a)\cup\sigma_{\Sc,\Ac}^{\Phi}(b)$. This proves that
Eq. (\ref{e0}).\\
The proof of $(\ref{e88})$ may be cheked in the same way as the proof of $(\ref{e0})$. To prove the inverse inclusion of $(\ref{e88})$. Let $q\not\in \sigma_{\Sc,\Ac}^{\Phi^{0}}(a+b)\cup\{0\}$. Since $ab$ and $ba$ are in $\Ac^{-1}(0)$, then by Eq. (\ref{e1}), we have
\begin{align*}R_{q}(a)R_{q}(b)=v_{0}+c_{0}\end{align*}
\noindent for some $v_{0}\in\Vc^{-1}$ and $c_{0}\in\Ac^{-1}(0)$. Also, since $\sigma_{\Sc,\Ac}^{\Phi^{0}}(a)=\sigma_{\Sc,\Ac}^{\Phi}(a)$, then there exists $v_{1}\in\Vc^{-1}$ and $c_{1}\in\Ac^{-1}(0)$ such that
\begin{align*}R_{q}(b)=v_{1}(c_{1}R_{q}(b)+v_{0}+c_{0}).\end{align*}
\noindent In this way, we see that $R_{q}(a)$ and $R_{q}(b)$ are in $\Phi_{\Ac}^{0}$. This proof is complete.\qed
\begin{theorem}\label{ts}
Let $a\in\Vc^{-1}$. Then
\begin{align*}\sigma_{\Sc,\Ac}^{\Phi}(a^{-1})=\left\{\displaystyle \frac{\overline{q}}{|q|^{2}}:\ q\in \sigma_{\Sc,\Ac}^{\Phi}(a)\right\}\end{align*}
 and
 \begin{align*}\sigma_{\Sc,\Ac}^{\Phi^{0}}(a^{-1})=\left\{\displaystyle \frac{\overline{q}}{|q|^{2}}:\ q\in \sigma_{\Sc,\Ac}^{\Phi^{0}}(a)\right\}.\end{align*}
\end{theorem}
\proof We can write for $q\in\mathbb{H}\backslash\{0\}$
\begin{align*}\displaystyle \frac{R_{q}(a)}{|q|^{2}}
&=\left(a^{-2}-\frac{2Re(q)}{|q|^{2}}a^{-1}+\frac{1}{|q|^{2}}\right)a^{2}\\
&=a^{2}\left(a^{-2}-\frac{2Re(q)}{|q|^{2}}a^{-1}+\frac{1}{|q|^{2}}\right)\end{align*}
because $\Vc$ is two-sided Banach algebra. This shows that $R_{q}(a)\in\Phi_{\Ac}$ if and only if $(a^{-2}-\frac{2Re(q)}{|q|^{2}}a^{-1}+\frac{1}{|q|^{2}})a^{2}\in\Phi_{\Ac}$. Since $a^{2}\in\Vc^{-1}$, then
\begin{align*}q\in \sigma_{\Sc,\Ac}^{\Phi}(a)\mbox{ if and only if }\frac{\overline{q}}{|q|^{2}}\in\sigma_{\Sc,\Ac}^{\Phi}(a^{-1}).\end{align*}
 Similary, we have that $R_{q}(a)\in\Phi_{\Ac}^{0}$ if and only if $R_{\frac{\overline{q}}{|q|^{2}}}\in \Phi_{\Ac}^{0}$.\qed
\begin{corollary}
Let $a,b\in\Vc^{-1}$.
\begin{enumerate}
\item If $a^{-1}-b^{-1}\in\Phi_{\Ac}$, then
\begin{align*}\sigma_{\Sc,\Ac}^{\Phi}(a)=\sigma_{\Sc,\Ac}^{\Phi}(b).\end{align*}
\item If $a^{-1}-b^{-1}\in\Ac^{-1}(0)$, then
\begin{align*}\sigma_{\Sc,\Ac}^{\Phi^{0}}(a)=\sigma_{\Sc,\Ac}^{\Phi^{0}}(b).\end{align*}
\end{enumerate}
\end{corollary}\label{perturbation}
\proof Since $a^{-1}-b^{-1}\in\Phi_{\Ac}$, then $\sigma_{\Sc,\Ac}^{\Phi}(a^{-1})=\sigma_{\Sc,\Ac}^{\Phi}(b^{-1})$. Applying Theorem \ref{ts},
we obtain $\sigma_{\Sc,\Ac}^{\Phi}(a)=\sigma_{\Sc,\Ac}^{\Phi}(b)$.\qed
\begin{theorem}
Let $\Ac: \Vc\longrightarrow\Wc$ be a bounded below homomorphism and $a\in\Vc$. Then,
\begin{align*}\partial\sigma_{S}(a)\subset\partial\sigma_{S,\Ac}^{\Phi}(a).\end{align*}
\end{theorem}
\proof Let $q\in \partial\sigma_{S}(a)$. Then $q\not\in int(\sigma_{S,\Ac}^{\Phi}(a))$. So, it suffices to prove that $q\in\sigma_{S,\Ac}^{\Phi}(a)$. Suppose that $R_{q}(a)\in\Phi_{\Ac}$. Then,
\begin{align*}(\Ac(a))^{2}-2Re(q)\Ac(a)+|q|^{2}1_{\Wc}\in\Wc^{-1}.\end{align*}
Now, let $(q_{n})_{n}$ be a sequence in $\mathbb{H}\backslash\sigma_{S}(a)$ such that $q_{n}\longrightarrow q$. So,
\begin{align*}(\Ac(a))^{2}-2Re(q_{n})\Ac(a)+|q_{n}|^{2}1_{\Wc}\longrightarrow(\Ac(a))^{2}-2Re(q)\Ac(a)+|q|^{2}1_{\Wc}\end{align*}
By the continuity of the map $u\longmapsto u^{-1}$, we get
\begin{align*}R_{q_{n}}^{-1}(\Ac(a))=\Ac(R_{q_{n}}^{-1}(a))\longrightarrow R_{q}^{-1}(\Ac(a)).\end{align*}
Therefore, the sequence $\Ac(R_{q_{n}}(a)))_{n}$ is Cauchy in $\Wc$. Since $\Ac$ is bounded below, then $(R_{q_{n}}^{-1}(a))_{n}$ converges, say $R_{q_{n}}^{-1}(a)\longrightarrow b\in\Vc$. This show that
\begin{align*}bR_{q}(a)=R_{q}(a)b=1_{\Vc}\end{align*}
which contradicts the fact that $q\in \sigma_{S}(a)$.\qed
\begin{lemma}\label{lem}
Let $(a_{n})$ be a sequence elements of $\Vc^{-1}$ converging to a non-invertible element. Then $\lim_{n}\|a_{n}^{-1}\|=+\infty$.
\end{lemma}
\proof The proof is the same as the proof in a complex Banach algebra, see \cite{AB}.\qed
\begin{theorem}
Let $a\in\Ac$ and $q\in \partial\sigma_{S}(a)$. Then there exists a sequence $(b_{n})$ of elements of $\Ac$ such that $\|b_{n}\|=1$ and
\begin{align*}\lim_{n}(a^{2}-2Re(q)a+|q|^{2}1_{\Vc})b_{n}=\lim_{n}b_{n}(a^{2}-2Re(q)a+|q|^{2}1_{\Vc})=0.\end{align*}
\end{theorem}
\proof Let $q\in\partial\sigma_{S}(a)$. Then, there exists a sequence $(q_{n})\subset\mathbb{H}\setminus\sigma_{S}(a)$ of quaternion number converging to $q$. Take
\begin{align*}b_{n}:=\displaystyle\frac{(a^{2}-2Re(q_{n})a+|q_{n}|^{2}1_{\Vc})^{-1}}{\|(a^{2}-2Re(q_{n})a+|q_{n}|^{2}1_{\Vc})^{-1}\|}.\end{align*}
Then we have $\|b_{n}\|=1$ and
\begin{align*}R_{q}(a)b_{n}=(a^{2}-2Re(q_{n})a+|q_{n}|^{2}1_{\Vc})b_{n}+(2Re(q_{n}-q)a+(|q^{2}|-|q_{n}|^{2})1_{\Vc})b_{n}.\end{align*}
Consequently
\begin{align*}\|R_{q}(a)b_{n}\|\leq\displaystyle\frac{1}{\|(a^{2}-2Re(q_{n})a+|q_{n}|^{2})^{-1}\|}+\|a\|\left(2|Re(q_{n}-q))|\right)+||q|^{2}-|q_{n}|^{2}|\end{align*}
By Lemma \ref{lem}, the result follow.\qed
\begin{theorem}
Let $\Ac:\Vc\longrightarrow\Wc$ be a homomorphism of Banach algebra and let $v_{1},v_{2}\in\Vc$. Take
\begin{align*}\mathbb{H}_{p,0}:=\left\{q\in\mathbb{H}^{*}:\ Re(q)=0\right\}.\end{align*}
Then,
\begin{align*}\sigma_{S,\Ac}^{\Phi}(v_{1}v_{2})\backslash\mathbb{H}_{p,0}=\sigma_{S,\Ac}^{\Phi}(v_{2}v_{1})\backslash\mathbb{H}_{p,0}.\end{align*}
\end{theorem}
\proof Let $a,b\in\Wc$ and $q\not\in\sigma_{S}(ab)\backslash\mathbb{H}_{p,0}$. Then there exist $c\in\Wc$ such that
\begin{align*}(abab+|q|^{2}1_{\Wc})c=c(abab+|q|^{2}1_{\Wc})=1_{\Wc}.\end{align*}
So, we have
\begin{align*}R_{q}(ba)(bcaba-1_{\Wc})
&=b(ab)^{2}caba-(ba)^{2}+|q|^{2}bcaba-|q|^{2}1_{\Wc}\\
&=b(1_{\Wc}-|q|^{2}c)aba-(ba)^{2}+|q|^{2}bcaba-|q|^{2}1_{\Wc}\\
&=-|q|^{2}1_{\Wc}\end{align*}
and
\begin{align*}(bcaba-1_{\Wc})R_{q}(ba)
&=bc(ab)^{2}aba+|q|^{2}bcaba-(ba)^{2}-|q|^{2}1_{\Ac}\\
&=-|q|^{2}1_{\Wc}\end{align*}
Consequently, $R_{q}(ba)$ is invertible in $\Wc$. Now, let $s=\Ac(v_{1})$ and $t=\Ac(v_{2})$, then

\begin{align*}\sigma_{S,\Ac}^{\Phi}(v_{1}v_{2})\backslash\mathbb{H}_{p,0}
&=\sigma_{S}^{\Phi}(\Ac(v_{1}v_{2}))\backslash\mathbb{H}_{p,0}\\
&=\sigma_{S}^{\Phi}(st)\backslash\mathbb{H}_{p,0}\\
&=\sigma_{S}^{\Phi}(ts)\backslash\mathbb{H}_{p,0}\\
&=\sigma_{S,\Ac}^{\Phi}(v_{2}v_{1})\backslash\mathbb{H}_{p,0}.\end{align*}\qed
\section{Boundary S-spectrum}
Let $\Vc$ be a quaternionic two-sided Banach algebra with unit. The boundary S-spectrum of $v\in\Vc$ is given by
\begin{align*}B_{S,\partial}(v):=\Big\{q\in\mathbb{H}:\ R_{q}(v)\in\partial(\Vc\backslash \Vc^{-1})\Big\}.\end{align*}
\begin{proposition}Let $\Vc$ be a quaternionic two-sided Banach algebra with unit. Then
\begin{align*}\partial\sigma_{S}(v)\subseteq B_{S,\partial}(v)\subseteq \sigma_{S}(v).\end{align*}
\end{proposition}
\proof Since $\Vc^{-1}\setminus \Vc$ is closed, $\partial(\Vc\backslash\Vc^{-1})\subseteq\Vc\backslash\Vc^{-1}$. So, $B_{S,\partial}(v)\subseteq\sigma_{S}(v)$. To prove $\partial\sigma_{S}(v)\subseteq B_{S,\partial}(v)$, let $q\in\partial\sigma_{S}(v)$. Then $R_{q}(v)\in\Vc\backslash\Vc^{-1}$. So, it suffice to prove that $R_{q}(v)\not\in int(\Vc\backslash\Vc^{-1})$. Since $q\not\in int(\sigma_{S}(v))$, then there exists a sequence $(q_{n})$ of elements of $\rho_{S}(v)$ such that $\lim_{n}q_{n}=q$. This show that
\begin{align*}\lim_{n}(v^{2}-2Re(q_{n})v+|q_{n}|)=R_{q}(v)\in\overline{\Vc^{-1}}\end{align*}
and therefore, $R_{q}(v)\not\in int(\Vc\backslash\Vc^{-1})$.\qed
\vskip 0.01 cm

In general, $\partial\sigma_{S}(v)$ is properly contained in $B_{S,\partial}(v)$, see Section \ref{shift}.
\begin{remark}\label{remarque1}Let $\Vc$ be a quaternionic two-sided Banach algebra with unit.
\begin{enumerate}
\item It follows, immediately, from the previous proposition that, for every $v\in\Vc$, the set $B_{S,\partial}(v)$ is non-empty compact set.
\item If $a\in\Vc^{-1}$, then $B_{S,\partial}(v^{-1})=(B_{S,\partial}(v))^{-1}$. Indeed, we can write for $q\neq 0$,
\begin{align*}R_{q}(v^{-1})=|q|^{2}R_{\frac{1}{q}}(v)v^{-2}=|q|^{2}v^{-2}R_{\frac{1}{q}}(v)\end{align*}
because $\Vc$ is two-sided Banach algebra. So, if $q\in B_{S,\partial}(v^{-1})$, then
\begin{align*}|q|^{2}R_{\frac{1}{q}}(v)v^{-2}\in \partial(\Vc\backslash\Vc^{-1}).\end{align*}
It follows from (2) that
\begin{align*}R_{\frac{1}{q}}(v)\in \partial(\Vc\backslash\Vc^{-1}).\end{align*}
\item Going over the same techniques of the proof of \cite[Theorem 2.15]{MS}, we obtain: if $\Ac:\ \Vc\longrightarrow\Wc$ is a continuous isomorphism then $\Ac(\partial(\Vc\setminus\Vc^{-1}))=\partial(\Wc\setminus\Wc^{-1})$.
\end{enumerate}
\end{remark}
\begin{theorem}
If $\Ac:\ \Vc\longrightarrow\Wc$ is a continuous isomorphism then
\begin{align*}B_{S,\partial}(v)=B_{S,\partial}(\Ac(v))=\displaystyle\cup_{c\in\Ac^{-1}(0)}B_{S,\partial}(v+c)\end{align*}
for all $v\in\Vc$.
\end{theorem}
\proof Let $q\in B_{S,\partial}(v)$, then $R_{q}(v)\in\partial(\Vc\setminus\Vc^{-1})$. therefore, by Remark \ref{remarque1} (4), we have
\begin{align*}(\Ac(v))^{2}-2Re(q)\Ac(v)+|q|^{2}\in\partial(\Wc\setminus\Wc^{-1}).\end{align*}
In this way we see that $q\in B_{S,\partial}(\Ac(v))$. Now, suppose that $q\in B_{S,\partial}(v+c)$ for some $c\in\Ac^{-1}(0)$, then
\begin{align*}R_{q}(\Ac(v))=\Ac((a+c)^{2}-2Re(q)(a+c)+|q|^{2})\in \partial(\Wc\setminus\Wc^{-1}).\end{align*}
Therefore, $q\in B_{S,\partial}(\Ac(v))$. Now, it remains to prove that $B_{S,\partial}(\Ac(v))\subset B_{S,\partial}(v)$. Let $q\in B_{S,\partial}(\Ac(v))$, then
\begin{align*}\Ac(R_{q}(v))=\Ac(v_{0})\mbox{ for some }v_{0}\in \partial(\Vc\setminus\Vc^{-1}).\end{align*}
 The injectivity of $\Ac$ imply that
 \begin{align*}v^{2}-2Re(q)v+|q|^{2}=v_{0}\in \partial(\Vc\setminus\Vc^{-1}).\end{align*}
 In this way we see that $q\in B_{S,\partial}(v)$.
\section{Application to right quaternionic linear operators}
\subsection{Preliminary}
In this subsection, we recall some definitions and we give some lemmas that we will need in the sequel. We refer to \cite{SL, RVA,KV} for more details.
\begin{definition}Let $\Vc_{\mathbb{H}}^{R}$ be a linear vector space under right multiplication by quaternionic scalars. For $f,g,h\in \Vc_{\mathbb{H}}^{R}$ and $q\in\mathbb{H}$, the inner product
\begin{align*}\langle .,.\rangle:\ \Vc_{\mathbb{H}}^{R}\times\Vc_{\mathbb{H}}^{R}\longrightarrow\mathbb{H}\end{align*}
 satisfies the following properties:
\begin{enumerate}
\item $\overline{\langle f,g\rangle}=\langle g,f\rangle$.
\item $\|f\|^{2}=\langle f,f\rangle>0$ unless $f=0$, a real norm.
\item $\langle f,g+h\rangle=\langle f,g\rangle+\langle f,h\rangle$.
\item $\langle f,gq\rangle=\langle f,g\rangle q$.
\item $\langle fq,g\rangle=\overline{q}\langle f,g\rangle$.
\end{enumerate}
\end{definition}
In the sequel, we assume that $\Vc_{\mathbb{H}}^{R}$ is complete under the norm given above and separable. In this case, $\langle .,.\rangle$ defines a right quaternionic Hilbert space.

\begin{proposition}\cite[Proposition 2.5]{RVA}
Let $\mathcal{F}=\{f_{k}:\ k\in\mathbb{N}\}$ be an orthonormal subset of $\Vc_{\mathbb{H}}^{R}$. Then, the following conditions are pairwise equivalent:
\begin{enumerate}
\item For every $f,g\in \Vc_{\mathbb{H}}^{R}$, the series $\sum_{k\in\mathbb{N}}\langle f,f_{k}\rangle\langle f_{k},g\rangle$ converges absolutely and it holds:
\begin{align*} \langle f,g\rangle=\sum_{k\in\mathbb{N}}\langle f,f_{k}\rangle\langle f_{k},g\rangle.\end{align*}
\item For every $f\in \Vc_{\mathbb{H}}^{R}$, it holds:
\begin{align*}\|f\|^{2}=\sum_{k\in\mathbb{N}}\mid \langle f_{k},f\rangle\mid^{2}\end{align*}
\item $\mathcal{F}^{\bot}=\{0\}$.
\item $\langle \mathcal{F}\rangle$ is dense in $\Vc_{\mathbb{H}}^{R}$.
\end{enumerate}
\end{proposition}
\begin{definition}
The set $\mathcal{F}$ as in the previous Proposition is called a Hilbert basis of $\Vc_{\mathbb{H}}^{R}$.
\end{definition}
\begin{proposition}\cite[Propsition 2.6]{RVA}
Every quaternionic Hilbert space $\Vc_{\mathbb{H}}^{R}$ has a Hilbert basis. All the Hilbert base of $\Vc_{\mathbb{H}}^{R}$ have the same cardinality.\\
Furthermore, if $\mathcal{F}$ is Hilbert basis of $\Vc_{\mathbb{H}}^{R}$, then every $f\in \Vc_{\mathbb{H}}^{R}$ can be uniquely decomposed as follows:
\begin{align*}f=\sum_{k\in\mathbb{N}}f_{k}\langle f_{k},f\rangle\end{align*}
where the series $\sum_{k\in\mathbb{N}}f_{k}\langle f_{k},f\rangle$ converges absolutely in $\Vc_{\mathbb{H}}^{R}$.
\end{proposition}
\begin{definition}
Let $\Vc_{\mathbb{H}}^{R}$ and $\mathcal{W}_{\mathbb{H}}^{R}$ be two quaternionic right Hilbert spaces. A mapping $T:\mathcal{D}(T)\subset\Vc_{\mathbb{H}}^{R}\longrightarrow\mathcal{W}_{\mathbb{H}}^{R}$, where $\mathcal{D}(T)$ stand for the domain of $T$, is called quaternionic right linear if
\begin{align*}T(fq+g)=T(f)q+T(g)\end{align*}
for all $f,g\in \Vc_{\mathbb{H}}^{R}$ and $q\in\mathbb{H}$.
\end{definition}
Denote by $\mathcal{L}(\Vc_{\mathbb{H}}^{R},\mathcal{W}_{\mathbb{H}}^{R})$ the set of all right linear operator from $\Vc_{\mathbb{H}}^{R}$ to $\mathcal{W}_{\mathbb{H}}^{R}$. If $\Vc_{\mathbb{H}}^{R}=\mathcal{W}_{\mathbb{H}}^{R}$ then $\mathcal{L}(\Vc_{\mathbb{H}}^{R},\Vc_{\mathbb{H}}^{R})$ is replaced by $\mathcal{L}(\Vc_{\mathbb{H}}^{R})$. We call an operator $T\in \mathcal{L}(\Vc_{\mathbb{H}}^{R},\mathcal{W}_{\mathbb{H}}^{R})$ bounded if there exist $M\in\mathbb{R}_{+}$ such that
\begin{align*}\|Tf\|\leq M\|f\|, \mbox{ for all }f\in\mathcal{D}(T).\end{align*}
Let $\mathcal{B}(\Vc_{\mathbb{H}}^{R},\mathcal{W}_{\mathbb{H}}^{R})$ denote the set of all bounded right linear operators from $\Vc_{\mathbb{H}}^{R}$ to $\mathcal{W}_{\mathbb{H}}^{R}$. If
 $\Vc_{\mathbb{H}}^{R}=\mathcal{W}_{\mathbb{H}}^{R}$, we will write  $\mathcal{B}(\Vc_{\mathbb{H}}^{R})$. The identity linear operator on $\Vc_{\mathbb{H}}^{R}$ will be denoted by $\mathbb{I}_{\Vc_{\mathbb{H}}^{R}}$. As in the complex case if $T\in \mathcal{L}(\Vc_{\mathbb{H}}^{R},\mathcal{W}_{\mathbb{H}}^{R})$, we set
 \begin{equation}\label{e9}\|T\|=\sup_{f\in\mathcal{D}(T)\backslash\{0\}}\displaystyle\frac{\|Tf\|}{\|f\|}.\end{equation}
 In this case, we have $T\in \mathcal{B}(\Vc_{\mathbb{H}}^{R},\mathcal{W}_{\mathbb{H}}^{R})$ if and only if $\|T\|<\infty$. Also, we have
 \begin{align*}\|T+S\|\leq\|T\|+\|S\|\mbox{ and }\|TS\|\leq\|T\|\|S\|.\end{align*}
 In the sequel, we investigate the $C^{*}-$algebras structure of $\mathcal{B}(\Vc_{\mathbb{H}}^{R})$. We refer to \cite{RVA} for more details. Set $\mathcal{F}=\{f_{k}:\ k\in\mathbb{N}\}$ the Hilbert basis of $\Vc_{\mathbb{H}}^{R}$. The left scalar multiplication on $\Vc_{\mathbb{H}}^{R}$ induced by $\mathcal{F}$ is defined as the map
 \begin{align*}
&\mathbb{H}\times \Vc_{\mathbb{H}}^{R}\longrightarrow \Vc_{\mathbb{H}}^{R}
\\
&\ (q,f)\longmapsto\ qf=\displaystyle\sum_{k\in\mathbb{N}}f_{k}q\langle f_{k},f\rangle.
\end{align*}
 The properties of the left scalar multiplication are described in the following proposition
 \begin{proposition}\cite[Proposition 3.1]{RVA} Let $f,g\in \Vc_{\mathbb{H}}^{R}$ and $p,q\in\mathbb{H}$, then
 \begin{enumerate}
 \item $q(f+g)=qf+qg$ and $q(fp)=(qf)p$.
 \item $\|qf\|=|q|\|f\|$.
 \item $q(pf)=(qp)f$.
 \item $\langle \overline{q}f,g\rangle=\langle f,qf\rangle$.
 \item $rf=fr$, for all $r\in\mathbb{R}$.
 \item $qf_{k}=f_{k}q$, for all $k\in\mathbb{N}$.
 \end{enumerate}
 It is easy to see that $(p+q)f=pf+qf$, for all $p,q\in\mathbb{H}$ and $f\in \Vc_{\mathbb{H}}^{R}$. We can also endow $\mathcal{B}(\Vc_{\mathbb{H}}^{R})$ with the two scalar multiplication
 \begin{equation}\label{e8}(qT)f=q(Tf)=\sum_{k\in\mathbb{N}}f_{k}q\langle f_{k},Tf\rangle\mbox{  and  }(Tq)f=T(qf),\end{equation}
 for all $q\in\mathbb{H}$ and $f\in \Vc_{\mathbb{H}}^{R}$.\\
 {\bf Notation }  Denote by $T^{*}$ the adjoint of $T\in\mathcal{B}(\Vc_{\mathbb{H}}^{R})$. We recall that the concept of adjoint quaternionic linear operator is the same as complex Hilbert space see \cite[Definition 2.12]{RVA}.
 \begin{theorem} \cite[Theorem 3.4]{RVA} Let $\Vc_{\mathbb{H}}^{R}$ be a right quaternionic Hilbert space equipped with a left scalar multiplication. Then the set $\mathcal{B}(\Vc_{\mathbb{H}}^{R})$, equipped pointwise sum, with the left and right scalar multiplication defined in Eq. (\ref{e8}), with the composition as product, with the adjunction $T\longrightarrow T^{*}$ as $*-$involution and with the norm defined in Eq. (\ref{e9}) is a quaternionic two-sided Banach $C^{*}-$algebras with unit $\mathbb{I}_{\Vc_{\mathbb{H}}^{R}}$.
 \end{theorem}
 The concept of quaternionic compact operators is the same as for complex Banach space. Let $\mathcal{K}(\Vc_{\mathbb{H}}^{R},\mathcal{W}_{\mathbb{H}}^{R})$
  denote the set of compact operators from $\Vc_{\mathbb{H}}^{R}$ to $\mathcal{W}_{\mathbb{H}}^{R}$. This class of operators was introduced and investigated in
   \cite{BK}. In particular, it is shown that $\mathcal{K}(\Vc_{\mathbb{H}}^{R})$ is closed two-sided ideal of $\mathcal{B}(\Vc_{\mathbb{H}}^{R})$, see \cite[Theorem 7.3]{BK}.
 \end{proposition}
 \subsection{Fredholm and Weyl S-spectra of quaternionic operators}
 We turn to the question of the S-spectrum. We will apply the results described in Section \ref{s:sum} to investigate the properties of the Fredholm and Weyl S-spectra of quaternionic operators. We recall that $\mathcal{K}(\Vc_{\mathbb{H}}^{R})$ is closed two-sided ideal of $\mathcal{B}(\Vc_{\mathbb{H}}^{R})$. So, $\mathcal{B}(\Vc_{\mathbb{H}}^{R})/\mathcal{K}(\Vc_{\mathbb{H}}^{R})$ is a unital Banach algebras with unit $\mathbb{I}_{\Vc_{\mathbb{H}}^{R}}+\mathcal{K}(\Vc_{\mathbb{H}}^{R})$. Now, we consider the natural quotion map
 \begin{align*}
&\pi: \mathcal{B}(\Vc_{\mathbb{H}}^{R}) \longrightarrow \mathcal{C}(\Vc_{\mathbb{H}}^{R}):=\mathcal{B}(\Vc_{\mathbb{H}}^{R})/\mathcal{K}(\Vc_{\mathbb{H}}^{R})
\\
&\ \ \quad \ \quad T \longmapsto\ \ [T]=T+\mathcal{K}(\Vc_{\mathbb{H}}^{R})
\end{align*}
 Note that $\pi$ is a unital homomorphism, see \cite{BK}. The norm on $\mathcal{C}(\Vc_{\mathbb{H}}^{R})$ is given by
 \begin{align*}\|[T]\|=\inf_{K\in\mathcal{K}(\Vc_{\mathbb{H}}^{R})}\|A+K\|.\end{align*}
 \begin{definition}
 The Calkin (or Fredholm) S-spectrum of $T\in\mathcal{B}(\Vc_{\mathbb{H}}^{R})$ is the  Fredholm S-spectrum of $T$ relative to the homomorphism $\pi$, i.e.
 \begin{align*}\sigma_{S,\pi}^{\Phi}(T)=\sigma_{S}(\pi(T))\end{align*}
 where $\sigma_{S}(\pi(T))$ is the S-spectrum of $\pi(T)$.
 \end{definition}
 Let $\sigma_{S,\pi}^{\Phi^{0}}(T)$ denote the Weyl S-spectrum of $T$ relative to the homorphism $\pi$. This two spectrum was introduced and studied by
   Muraleetharam and  Thirulogasanthar in \cite{BK,BK2}. In particular, we have $A\in \Phi_{\pi}^{0}$ if and only if $A\in\Phi_{\pi}$ and $ind(A)=0$, where
 \begin{align*}ind(A)=dim(A^{-1}(0))-dim((A^{*})^{-1}(0)).\end{align*}
 \vskip 0.1 cm

 \begin{theorem}\label{sum1}
Let $A,\ B\in\mathcal{B}(\Vc_{\mathbb{H}}^{R})$.
\begin{enumerate}
\item If $AB\in \mathcal{K}(\Vc_{\mathbb{H}}^{R})$ and $BA\in \mathcal{K}(\Vc_{\mathbb{H}}^{R})$, then
\begin{equation}\sigma_{\Sc,\pi}^{\Phi}(A+B)\backslash\{0\}=
\left[\sigma_{\Sc,\pi}^{\Phi}(A)\cup\sigma_{\Sc,\pi}^{\Phi}(B)\right]\backslash\{0\}.\end{equation}
and
\begin{equation}\sigma_{\Sc,\pi}^{\Phi^{0}}(A+B)\backslash\{0\}\subset
\left[\sigma_{\Sc,\pi}^{\Phi^{0}}(A)\cup\sigma_{\Sc,\pi}^{\Phi^{0}}(B)\right]\backslash\{0\}.\end{equation}
\noindent If, further, $\sigma_{\Sc,\pi}^{\Phi^{0}}(A)=\sigma_{\Sc,\pi}^{\Phi}(A)$, then
\begin{equation}\sigma_{\Sc,\pi}^{\Phi^{0}}(A+B)\backslash\{0\}=\left[\sigma_{\Sc,\pi}^{\Phi^{0}}(A)
\cup\sigma_{\Sc,\Ac}^{\Phi^{0}}(B)\right]\backslash\{0\}.\end{equation}
\item If $0\not\in\sigma_{S}(A)\cup\sigma_{S}(B)$ and $A^{-1}-B^{-1}\in \mathcal{K}(\Vc_{\mathbb{H}}^{R})$, then
\begin{align*}\sigma_{S,\pi}^{\Phi}(A)=\sigma_{S,\pi}^{\Phi}(B)\mbox{  and  }\sigma_{S,\pi}^{\Phi^{0}}(A)=\sigma_{S,\pi}^{\Phi^{0}}(B).\end{align*}
\item Recall that $\mathbb{H}=\{q\in\mathbb{H}^{*}:\ Re(q)=0\}$. Then
\begin{align*}\sigma_{S,\pi}^{\Phi}(AB)\backslash\mathbb{H}_{p,0}=\sigma_{S,\pi}^{\Phi}(BA)\backslash\mathbb{H}_{p,0}.\end{align*}
\end{enumerate}
\end{theorem}
\proof Note that
   \begin{align*}\mathcal{K}(\Vc_{\mathbb{H}}^{R})\subset Pr(\Phi_{\pi})\mbox{   and   }\mathcal{K}(\Vc_{\mathbb{H}}^{R})=\pi^{-1}([0]).\end{align*}
 Applying Theorem \ref{sum} and Corollary \ref{perturbation} the results follows.\qed
 \vskip 0.2 cm

We mention two remarks.
\begin{remark}
For the complex Banach space, in \cite {BJ}, we have prove that if $T_{1},T_{2}\in \mathcal{B}(X)$ such that $\sigma_{ess}(T_{1}T_{2})=\{0\}$ and $T_{1}T_{2}-T_{2}T_{1}\in\mathcal{K}(X)$, then
\begin{align*}\sigma_{ess}(T_{1}+T_{2})\backslash\{0\}=
\left[\sigma_{ess}(T_{1})\cup\sigma_{ess}(T_{2})\right]\backslash\{0\},\end{align*}
where $\sigma_{ess}(T_{i})$ is the essential spectrum of $T_{i}$. The case of quaternionic Hilbert space is more complicated and we were not able to refine the hypothesis in the previous theorem. We leave this question open.
\end{remark}
\begin{remark}
 If $K\in\mathcal{K}(\Vc_{\mathbb{H}}^{R})$ then $\sigma_{\Sc,\pi}^{\Phi^{0}}(K)=\sigma_{\Sc,\pi}^{\Phi}(K)=\{0\}$.
\end{remark}
In the wat follows we show that there exist a quaternionic operator such that its Calkin S-spectrum coincides with its  Weyl S-spectrum. We first mention the following theorem:
\begin{theorem}\label{st}\cite[Theorem 6.20]{BK}
Let $T\in\Phi_{\pi}$. Then there exists a constant $\varepsilon>0$ such that for every operators $S\in\mathcal{B}(\Vc_{\mathbb{H}}^{R})$ with
$\|S\|<\varepsilon,\ A+S\in \Phi_{\pi}$ and $ind(A+S)=ind(A)$.
\end{theorem}
Set:
\begin{align*}\Phi_{T,\pi}:=\Big\{q\in\mathbb{H}:\ R_{q}(T)\in \Phi_{\pi}\Big\}.\end{align*}
\begin{theorem}\label{res1}Let $T$ be a bounded right operator. Then, $q\longmapsto ind(R_{q}(T))$
is constant on any component of $\Phi_{T,\pi}$.\end{theorem}
\proof Let $q_{1}$ and $q_{2}$ be two elements of $\Phi_{T,\pi}$ that are connected by a smooth curve $\Gamma$. Take $q\in \Gamma$, then $R_{q}(T)\in \Phi_{\pi}$. By Theorem \ref{st}, there exists a constant $\varepsilon>0$ such that for all $S\in\mathcal{B}(\Vc_{\mathbb{H}}^{R})$ with
$\|S\|<\varepsilon$,
\begin{align*}R_{q}(T)+S\in\Phi_{\pi}\mbox{  and  }ind(R_{q}(T)+S)=ind(R_{q}(T)).\end{align*}
Set
\begin{align*}\Oc(A,q,\varepsilon):=\Big\{q'\in\mathbb{H}:\ 2\mid Re(q)-Re(q')\mid\|A\|+\mid |q'|^{2}-|q|^{2}\mid<\varepsilon.\Big\}\end{align*}
Let $q'\in \Oc(A,q,\varepsilon)$. Then
\begin{align*}R_{q}(A)+2(Re(q)-Re(q'))A+ (|q'|^{2}-|q|^{2})\mathbb{I}_{\Vc_{\mathbb{H}}^{R}}=R_{q'}(A)\in\Phi_{\pi}\end{align*}
and
\begin{align*}ind(R_{q}(A)+2(Re(q)-Re(q'))A+ (|q'|^{2}-|q|^{2})\mathbb{I}_{\Vc_{\mathbb{H}}^{R}})=ind(R_{q}(A)).\end{align*}
\noindent A cover of the set $\Gamma$ by open sets was thus constructed. By the Heine-Borel theorem, there is a finite number of such sets which cover $\Gamma$. Since the index $q\longmapsto ind(R_{q}(A))$ is constant on each of these sets (each crosses at least on other), then $ind(R_{q_{1}}(A))=ind(R_{q_{2}}(A))$.
\begin{corollary}\label{res2}
Let $T$ be a bounded right operator. If $\Phi_{T,\pi}$ is connected, then
\begin{align*}\sigma_{S,\pi}^{\Phi}(T)=\sigma_{S,\pi}^{\Phi^{0}}(T).\end{align*}
\end{corollary}
\proof Since the inclusion $\sigma_{S,\pi}^{\Phi}(T)\subset\sigma_{S,\pi}^{\Phi^{0}}(T)$ is know, then it suffices to show that $\Phi_{T,\pi}\cap\sigma_{S,\pi}^{\Phi^{0}}(T)$ is empty. Let $q\in\Phi_{T,\pi}$. Since $\rho_{S}(T)$ is not empty, then there exists $p\in\mathbb{H}$ such that $R_{p}(T)\in \mathcal{B}(\Vc_{\mathbb{H}}^{R})^{-1}$ and consequently $R_{p}(T)\in \Phi_{\pi}$ and $ind(R_{p}(T))=0$. Moreover, $\Phi_{T,\pi}$ is connected, by the previous theorem, $ind(T)$ is constant on $\Phi_{T,\pi}$. Therefore
 \begin{align*}ind(R_{q}(T))=ind(R_{q}(T))=0.\end{align*}
In this way we see that $q\not\in \sigma_{S,\pi}^{\Phi^{0}}(T)$.\qed
\begin{corollary}
Let $T\in \mathcal{B}(\Vc_{\mathbb{H}}^{R})$ and $S\in \mathcal{B}(\Vc_{\mathbb{H}}^{R})$. Assume that $ST$ and $TS$ are right linear compact operators, then
\begin{equation*}\label{e99}\sigma_{\Sc,\pi}^{\Phi^{0}}(T+S)\backslash\{0\}=\left[\sigma_{\Sc,\pi}^{\Phi^{0}}(T)
\cup\sigma_{\Sc,\Ac}^{\Phi^{0}}(S)\right]\backslash\{0\}.\end{equation*}
\end{corollary}
\proof Combine Theorem \ref{res1} and Corollary \ref{res2}.\qed
\begin{corollary}
Let $T$ be a bounded right operator. If $\Phi_{T,\pi}$ is connected, then
\begin{align*}\sigma_{k}^{S}(T)=\emptyset\end{align*}
for all $k\in\mathbb{Z}\backslash\{0\}$, where
\begin{align*}\sigma_{k}(T)=\{q\in\mathbb{H}:\ R_{q}(T)\in\Phi_{\pi}\mbox{ and }ind(R_{q}(T))=k\}.\end{align*}
\end{corollary}
\proof Combine Corollary \ref{res2} with \cite[Remark 6.7]{BK2}.\qed

\subsection{Right shift}\label{shift}
We consider the right quaternionic space:
\begin{align*}\ell^{2}_{\mathbb{H}}(\mathbb{Z})=\Big\{x:\mathbb{Z}\longrightarrow\mathbb{H}\mbox{  such that  } \|x\|^{2}:=\sum_{i\in\mathbb{Z}}|x_{i}|^{2}<\infty\Big\}.\end{align*}
with the right left scalar multiplication
\begin{align*}xa=(x_{i}a)_{i\in\mathbb{Z}}\end{align*}
for $x=(x_{i})_{i\in\mathbb{Z}}$ and $a\in\mathbb{H}$. The associated scalar product is given by
\begin{align*}\langle x,y\rangle:=\langle x,y\rangle_{\ell^{2}_{\mathbb{H}}(\mathbb{Z})}:=\sum_{i\in\mathbb{Z}}\overline{x_{i}}y_{i}.\end{align*}
Take $\{e_{n}:\ n\in\mathbb{Z}\}$ the Hilbert basis of $\ell^{2}_{\mathbb{H}}(\mathbb{Z})$ where, for each $n\in\mathbb{Z}$, $e_{n}=(e_{n}^{k})_{k\in\mathbb{Z}}$, where $e_{n}^{n}=1$ and $e_{n}^{i}=0$ for all $i\neq n$. Consider the right shift
 \begin{align*}
R: &\ell_{\mathbb{H}}^{2}(\mathbb{Z})\longrightarrow \ell_{\mathbb{H}}^{2}(\mathbb{Z})
\\
& x\longmapsto\ y=(y_{i})_{i\in\mathbb{Z}}
\end{align*}
where $y_{i}=x_{i+1}$ if $i\neq -1$ and $0$ if $i=-1$. We have $\|R(x)\|^{2}=\sum_{i\neq -1}|x_{i}|^{2}\leq \|x\|^{2}$. Also, since $\|R(e_{1})\|=1$, we have $\|R\|=1$. We recall the following definition:
\begin{definition}\cite{BK}
Let $A\in \mathcal{B}(V_{\mathbb{H}}^{R})$. The approximate $S-$spectrum of $A$, denoted by $\sigma_{app}^{S}(A)$, is defined as:
\begin{align*}\sigma_{app}^{S}(A)=\Big\{q\in\mathbb{H}:
\exists (f_{n})_{n} \mbox{ with } \|f_{n}\|=1\mbox{  and  }\|R_{q}(A)f_{n}\|\rightarrow 0\Big\}.\end{align*}
\end{definition}
\vskip 0.1 cm

The set of right eigenvalues coincide with the point $S-$spectrum, see \cite{RVA}, Proposition 4.5. We turn to the question of the spectrum of the right shift $R$. Consider the right eigenvalue problem  of $R$. Let $q\in\mathbb{H}$ and $x=(x_{i})_{i\in\mathbb{Z}}\in\ell^{2}_{\mathbb{H}}(\mathbb{Z})$ with $x_{-i}=0$ for all $i\in\mathbb{N}^{*}$ such that $Rx=xq$, then $x_{i+1}=x_{i}q$ for all $i\in\mathbb{N}$. This implies that $B_{\mathbb{H}}(0,1)\subset\sigma_{app}^{S}(R)\subset\sigma_{S}(R)$. Going over same technique of \cite[Example 6.31]{BK4} we have $\sigma_{S}(R)=\nabla_{\mathbb{H}}(0,1)$ the closed quaternionic unit ball. In particular $0\not\in \partial\sigma_{S}(R)$: the boundary of $\sigma_{S}(A)$. In the following, we will show that  $0\in B_{S,\partial}(R)$ the  boundary S-spectrum of $R$. To begin, let's consider the map:
\begin{align*}
T: &\ell_{\mathbb{H}}^{2}(\mathbb{Z})\longrightarrow \ell_{\mathbb{H}}^{2}(\mathbb{Z})
\\
& x\longmapsto\ y=(y_{i})_{i\in\mathbb{Z}}
\end{align*}
where $y_{-1}=x_{0}$ and $y_{j}=0$ for all  $j\neq -1$. We note that $T(e_{0})=e_{-1}$ and $T(e_{k})=0$ for all $k\neq 0$. Therefore,
\begin{align*}(Tq)(x)
&=T(\sum_{i\in\mathbb{Z}}e_{k}q\langle e_{k},x\rangle)\\
&=T(e_{0}q\langle e_{0},x \rangle)\\
&=e_{-1}qx_{0}\end{align*}
\noindent for all $q\in\mathbb{H}$ and $x=(x_{i})_{i\in\mathbb{Z}}\in\ell^{2}_{\mathbb{H}}(\mathbb{Z})$.
 \begin{proposition}
 $\sigma_{S,\pi}^{\Phi^{0}}(R)$ is properly contained in $\sigma_{S}(R)$.
 \end{proposition}
 \proof We have $(R+T)(x)=y=(y_{i})_{i\in\mathbb{Z}}$, where $y_{i}=x_{i+1}$ for all $i\in\mathbb{Z}$. Therefore,
  \begin{align*}R+T\in\Phi_{\pi}\mbox{ and }ind(R+T)=0.\end{align*}
   Since $T\in\mathcal{K}(\ell_{\mathbb{H}}^{2}(\mathbb{Z}))$, we can apply \cite[Theorem 6.16]{BK}, we infer that
   \begin{align*}R\in\Phi_{\pi}\mbox{  and  } ind(R)=ind(R+T)=0.\end{align*}
 Again, using \cite[Corllary 6.15]{BK}, we have
 \begin{align*}R^{2}\in\Phi_{\pi}\mbox{ and }ind(R^{2})=2ind(R)=0.\end{align*}
 Hence $0\in  \sigma_{S}(R)\backslash\sigma_{S,\pi}^{\Phi^{0}}(R)$.\qed

 The technique of the proof of the following Proposition are inspired from \cite{MS} in the case of $\ell_{\mathbb{C}}^{2}(\mathbb{Z})$.
\begin{proposition}
$\partial\sigma_{S}(R)$ is properly contained in $B_{S,\partial}(R)$.
\end{proposition}
\proof For $q\in\mathbb{H}\backslash\{0\}$, $(R+Tq)^{2}\in(\mathcal{B}(\ell_{\mathbb{H}}^{2}(\mathbb{Z})))^{-1}$ and its inverse is given by
\begin{align*}(R+Tq)^{-2}y=x=(x_{i})_{i\in\mathbb{Z}},\end{align*}
where
\begin{align*}x_{i}=\displaystyle\left\{ \begin{array}{rl}
\noindent\frac{q}{|q|^{2}}y_{i-2},  \mbox{ if } i\in\{0,1\},
\\~~
\\
 y_{i-2},\ \ \ \ \ \
  \mbox{ otherwise }.
\end{array} \right.\end{align*}
\noindent  Now, let $0<|q|<1$, then $R_{n}=(R+Tq^{n})^{2}$ is invertible for all $n\in\mathbb{N}$ and
\begin{align*}\| R_{n}-R^{2}\|=\|RTq^{n}+Tq^{n}R+(Tq^{n})^{2}\|\leq 2|q|^{n}+|q|^{2n}.\end{align*}
So, $R_{n}$ converges to $R^{2}$. This imply that  $R^{2}\in \mathcal{B}(\ell_{\mathbb{H}}^{2}(\mathbb{Z}))\backslash\overline{\mathcal{B}^{-1}(\ell_{\mathbb{H}}^{2}(\mathbb{Z}))}$.\quad

\end{document}